\documentclass[11pt]{article}
\usepackage{amsfonts}
\usepackage{mathrsfs}
\usepackage{amsmath}
\usepackage{amssymb}
\usepackage{graphicx}
\usepackage{epsfig}
\usepackage{epic}
\usepackage{cite}
\usepackage{setspace}
\usepackage{amsthm,latexsym}
\usepackage{float}

\renewcommand{\paragraph}{\roman{paragraph}}
\newtheorem{theorem}{\scshape \mdseries \bf Theorem}[section]
\newtheorem{lemma}[theorem]{\scshape \mdseries  \bf Lemma}

\newtheorem{corol}[theorem]{\scshape \mdseries  \bf Corollary}

\begin{document}

\title{\sf Automorphisms of the subspace sum graphs on a vector space}
\author{ Fenglei Tian, \ \  Dein Wong \thanks{Corresponding
author. \ \ E-mail address: wongdein@163.com.
Supported by ``the Fundamental Research Funds for the Central Universities
(No. 2017BSCXB53)''. }
\\ {\small  \it  School of Mathematics, China University of Mining and
Technology, Xuzhou 221116, China.}
   }
\date{}
\maketitle

\noindent {\bf Abstract:}\ The subspace sum graph $\mathcal{G}(\mathbb{V})$ on a finite dimensional vector space $\mathbb{V}$ was introduced by Das [Subspace Sum Graph of a Vector Space, arXiv:1702.08245], recently. The vertex set of $\mathcal{G}(\mathbb{V})$ consists of all the nontrivial proper subspaces of $\mathbb{V}$ and two distinct vertices $W_1$ and $W_2$ are adjacent if and only if $W_1+W_2=\mathbb{V}$. In that paper, some structural indices (e.g., diameter, girth, connectivity, domination number, clique number and chromatic number) were studied, but the characterization of automorphisms of $\mathcal{G}(\mathbb{V})$ was left as one of further research topics. Motivated by this, we in this paper characterize the automorphisms of $\mathcal{G}(\mathbb{V})$ completely.



\section{Introduction}
\quad
Up to now, there are some abstract graphs defined on certain algebraic objects, which are proposed from the motivation that the properties of those algebraic structures can be revealed by studying the graphs associated with them. For example, the zero-divisor graphs defined on finite rings \cite{Akbari3}, commutative rings \cite{Anderson2} or noncommutative rings \cite{Akbari2}, regular graphs or total graphs on a commutative ring \cite{Akbari1, Anderson1}, cozero-divisor graph of a commutative ring \cite{Afkhami} and so on. For various graphs defined on finite dimensional vector space, readers can refer to \cite{Das2,Das3,Das4,Das5} for details.

Let $\mathbb{V}$ be a finite dimensional vector space of dimension greater than one over a field $F$. Recently, the subspace sum graph $\mathcal{G}(\mathbb{V})$ on $\mathbb{V}$ was introduced by Das \cite{Das1}, the vertex set of which is composed of all the nontrivial proper subspaces of $\mathbb{V}$ and two vertices $W_1, W_2$ are adjacent, written as $W_1\sim W_2$, if and only if $W_1+W_2=\mathbb{V}$. In \cite{Das1}, the author investigated the diameter, girth, connectivity, domination number, clique number, chromatic number of $\mathcal{G}(\mathbb{V})$ and the properties of $\mathcal{G}(\mathbb{V})$ with the base field $F$ being finite.

Since automorphisms of graphs can reveal the relationship among their vertices, then they can help analyze the structure of graphs and let us go further to achieve the motivation stated at the beginning. Moreover, automorphisms of graphs are also of importance in algebraic graph theory. Thus, characterizing the automorphisms of graphs also have attracted many attentions.
Let $R$ be a ring and $Z(R)\setminus\{0\}$ the nonzero zero divisor set of $R$. The zero-divisor graph $\Gamma (R)$ of $R$ is defined as a graph with vertex set $Z(R)\setminus\{0\}$ and for vertices $x, y$, there is a directed edge from $x$ to $y$ if $xy=0$. Anderson and Livingston \cite{Anderson2} proved that the automorphism group of $\Gamma(\mathbb{Z}_n)$ ($n\geq 4$) is a direct product of some symmetric groups. Wong et al. \cite{Wong} characterized the automorphisms of the zero-divisor graph, whose vertex set consists of all rank one upper triangular matrices over a finite field. Applying the main results of \cite{Wong}, Wang \cite{Wang1,Wang2} obtained the automorphisms of the zero-divisor graph on upper triangular matrix ring or  full matrix ring with a finite base field, respectively.
Inspired by these results, we can ask a natural problem: {\sf how about the automorphisms of a subspace sum graph?} In addition, this is also an open problem for further research on subspace sum graphs in \cite{Das1}.

Hence, in present paper we focus on this problem and address it completely (see Section 3).  Prior to presenting the proof of it, some lemmas are demonstrated in Section 2.

\section{Preliminaries and Lemmas}
\quad
Let $F_q$ be a finite field of order $q$ and $q=p^m$ with $p$ a prime integer.
Throughout, we let $\mathbb{V}$ be a finite dimensional vector space with the base field $F_q$. The dimension of a subspace $W$ of $\mathbb{V}$ is denoted by $dim(W)$. Denote by $V$ the vertex set of $\mathcal{G}(\mathbb{V})$. From the definition of graph $\mathcal{G}(\mathbb{V})$, it is clear that if $dim(\mathbb{V})=1$, then the vertex set $V$ is empty. If $dim(\mathbb{V})=2$, one can easily see that $\mathcal{G}(\mathbb{V})$ is a complete graph (also see Theorem 3.1 of \cite{Das1}). Hence, in the following, we always set $dim(\mathbb{V})=n\geq 3$. Let $e_i\in \mathbb{V}$ be a vector with the $i$-th component 1 and the others 0. Then $\{e_1, e_2, \ldots, e_n\}$ constructs a basis of $\mathbb{V}$ and any vector $\alpha\in \mathbb{V}$ can be uniquely expressed as $\alpha=\sum_{i=1}^{n}a_ie_i$ with $a_i\in F_q$.

First, we present two standard automorphisms of $\mathcal{G}(\mathbb{V})$ as follows.

{\bf Invertible linear transformation}

Let $X=[x_{ij}]$ be an invertible matrix of order $n$ over $F_q$. Define a mapping $\sigma_X$ from $\mathbb{V}$ to itself as $$\sigma_X(\alpha)=\sigma_X(\sum_{i=1}^{n}a_ie_i)=\sum_{i=1}^{n}(\sum_{j=1}^{n}x_{ij}a_j)e_i, \forall \alpha=\sum_{i=1}^{n}a_ie_i\in \mathbb{V}.$$
Clearly, $\sigma_X$ is an invertible linear transformation of $\mathbb{V}$, which naturally generates an automorphism of $\mathcal{G}(\mathbb{V})$, also written as $\sigma_X$, such that $\sigma_X(W)=\{\sigma_X(w), w\in W\}$ for $W\in V$.

{\bf Field automorphism}

Let $f$ be an automorphism of the base field $F_q$ and $\sigma_f$ a mapping from $\mathbb{V}$ to itself such that
$$\sigma_f(\beta)=\sigma_f(\sum_{i=1}^{n}b_ie_i)=\sum_{i=1}^{n}\sigma_f(b_i)e_i, \forall \beta=\sum_{i=1}^{n}b_ie_i \in \mathbb{V}.$$
Then one can easily check that the mapping on $V$ generated by $\sigma_f$, also denoted by $\sigma_f$, sending $W\in V$ to $\{\sigma_f(w), w\in W\}$ is an automorphism of $\mathcal{G}(\mathbb{V})$.

Next, we give some lemmas which will be used later.

\begin{lemma}{\rm (Theorem 6.1, \cite{Das1})} \ Let $W$ be a $k$-dimensional subspace of $\mathbb{V}$ over the base field $F_q$. Then the degree of $W$ in $\mathcal{G}(\mathbb{V})$, denoted by $deg(W)$, is $deg(W)=\sum_{r=0}^{k-1}N_r$, where
$$N_r=\frac{(q^k-1)(q^k-q)\cdots (q^k-q^{r-1})(q^n-q^k)(q^n-q^{k+1})\cdots (q^n-q^{n-1})}{(q^{n-k+r}-1)(q^{n-k+r}-q)\cdots (q^{n-k+r}-q^{n-k+r-1})}.$$
\end{lemma}

From Lemma 2.1, it is clear that $deg(W)$ only depends on the dimension $k$ of $W$, and if $W_1, W_2\in V$ are of distinct dimensions $k_1, k_2$, then $deg(W_1)\neq deg(W_2)$. Thus Lemma 2.2 follows from Lemma 2.1 immediately.

\begin{lemma}\ Let $\sigma$ be an automorphism of $\mathcal{G}(\mathbb{V})$. Then $\sigma$ sends each $k$-dimensional subspace of $\mathbb{V}$ to a subspace of equal dimension, where $1\leq k\leq n-1$.
\end{lemma}

Let $S$ be a subset of $\mathbb{V}$ and the subspace spanned by $S$ be $\langle S \rangle$. The $1$-dimensional subspace spanned by a nonzero vector $\alpha$ is denoted by $\langle \alpha \rangle$.

\begin{lemma}\ Let $\sigma$ be an automorphism of $\mathcal{G}(\mathbb{V})$ and $W$ a $k$-dimensional subspace ($2\leq k\leq n-1$). Suppose $\{\alpha_1, \alpha_2, \ldots, \alpha_k\}$ is a basis of $W$ and $\sigma(\langle \alpha_i \rangle)=\langle \beta_i \rangle$ for $1\leq i\leq k$. Then we have the following conclusions.
\end{lemma}
\begin{spacing}{1}
\begin{enumerate}
\item[(i)]  \textit{$\langle \beta_i \rangle\subset \sigma(W)$ for $1\leq i\leq k$};
\item[(ii)]  \textit{$\beta_1, \beta_2, \ldots, \beta_k$ are linearly independent};
\item[(iii)]  \textit{$\sigma(W)=\langle\beta_1, \ldots, \beta_k \rangle$}.
\end{enumerate}
\end{spacing}

\noindent
{\bf Proof.} \ First, we show the proof of (i). Assume on the contrary that there exists $\langle \beta_i \rangle$ which is not a subspace of $\sigma(W)$. From Lemma 2.2, we see $dim(\sigma(W))=k$, then let $\{\gamma_1, \gamma_2, \ldots, \gamma_k\}$ be a basis of $\sigma(W)$. From the assumption, we obtain that $\beta_i, \gamma_1, \gamma_2, \ldots, \gamma_k$ are linearly independent. Extend $\{\beta_i, \gamma_1, \gamma_2, \ldots, \gamma_k\}$ to a basis $\{\beta_i, \gamma_1, \gamma_2, \ldots, \gamma_k, \ldots, \gamma_{n-1}\}$ of $\mathbb{V}$. Let $W_1=\langle \gamma_1, \gamma_2, \ldots, \gamma_k, \ldots, \gamma_{n-1} \rangle$ be a $(k-1)$-dimensional subspace, then $\sigma(W)\subset W_1$ and $\sigma(W)+W_1=W_1\neq \mathbb{V}$. Thus $\sigma(W)$ is not adjacent to $W_1$, i.e., $\sigma(W)\nsim W_1$. Furthermore, $\langle \beta_i \rangle+ W_1=\mathbb{V}$, and thus $\langle \beta_i \rangle \sim W_1$. Since $\sigma$ is an automorphism, so is $\sigma^{-1}$. Then applying $\sigma^{-1}$ to $\sigma(W)\nsim W_1$ and $\langle \beta_i \rangle \sim W_1$, we get $W\nsim \sigma^{-1}(W_1)$ and $\langle \alpha_i \rangle \sim \sigma^{-1}(W_1)$, which contradicts with $\langle \alpha_i \rangle\subset W$. As a result, the assumption does not hold, and the conclusion follows.

\vskip 3pt
Second, We will apply induction on $k$ to complete the proof of (ii). First, if $k=2$ and $\beta_1, \beta_2$ are not linearly independent, then $\beta_1=a\beta_2$ with $0\neq a\in F_q$. Thus $\langle \beta_1 \rangle=\langle \beta_2 \rangle$, that is, $\sigma(\langle \alpha_1 \rangle)=\sigma(\langle \alpha_2 \rangle)$, a contradiction.

Now suppose the conclusion holds for any $k=s-1\geq 2$ and we will show the case when $k=s$.
Let $\alpha_1, \alpha_2, \ldots, \alpha_s$ ($s> 2$) are linearly independent vectors and $\sigma(\langle \alpha_i \rangle)=\langle \beta_i \rangle$ for $1\leq i\leq s$. Extend $\{\alpha_1, \alpha_2, \ldots, \alpha_s\}$ to a basis $\{\alpha_1, \alpha_2, \ldots, \alpha_s, \alpha_{s+1}, \ldots, \alpha_n\}$ of $\mathbb{V}$. Let $W_2=\langle \alpha_s \rangle$ and $W_3=\langle \alpha_1, \ldots, \alpha_{s-1}, \alpha_{s+1}, \ldots, \alpha_n \rangle$, then
$$W_2+W_3=\mathbb{V}\ (i.e., W_2\sim W_3).$$
From the induction hypothesis, we know that $\beta_1, \beta_2, \ldots, \beta_{s-1}$ are linearly independent. If $\beta_1, \ldots, \beta_{s-1}, \beta_s$ are not linearly independent, then $\beta_s$ can be linearly expressed by $\beta_1, \ldots, \beta_{s-1}$, that is, $\beta_s= \sum_{i=1}^{s-1}a_i\beta_i$ with $a_i\in F_q$ and $\sum_i |a_i|\neq 0$.  Then
$$\sigma(W_2)=\langle \beta_s \rangle \subset \langle \beta_1, \ldots, \beta_{s-1}\rangle.$$
In addition, by the result of (i), $\langle \beta_i \rangle\subset \sigma(W_3)$ for $1\leq i\leq s-1$.
Hence $$\langle \beta_1, \ldots, \beta_{s-1}\rangle \subset \sigma(W_3).$$ Consequently, we derive that
$$\sigma(W_2)\subset \sigma(W_3),$$ which implies that $W_2\nsim W_3$, a contradiction. So, we say $\beta_1, \ldots, \beta_{s-1}, \beta_s$ are linearly independent.

\vskip 3pt
At last, we prove the result of (iii). Combining (i) and (ii), we see that $\langle \beta_1, \ldots, \beta_k\rangle \subset \sigma(W)$ and $dim(\langle \beta_1, \ldots, \beta_k\rangle)=k=dim(W)$, which complete the proof.\hfill$\square$

\vskip 3pt
The following lemma is a more general result of Lemma 2.3 (i).

\begin{lemma}\ Let $\sigma$ be an automorphism of $\mathcal{G}(\mathbb{V})$ and $W=\langle \alpha_1, \alpha_2, \ldots, \alpha_k \rangle$ ($2\leq k\leq n-1$) a $k$-dimensional subspace. Suppose $\alpha$ is a nonzero vector, we obtain
\end{lemma}
\begin{spacing}{1}
\begin{enumerate}
\item[(i)]  \textit{if $\langle \alpha \rangle \subset W$, then $\sigma(\langle \alpha \rangle) \subset \sigma(W)$};
\item[(ii)]  \textit{if $\langle \alpha \rangle$ is not a subspace of $W$, then $\sigma(\langle \alpha \rangle)$ is not a subspace of $\sigma(W)$}.
\end{enumerate}
\end{spacing}

\noindent
{\bf Proof.} \ Here we mainly give the proof of (i) and the other can be proved similarly.
Since $\langle \alpha \rangle \subset W$, then $\alpha$ can be linearly expressed by $\alpha_1, \alpha_2, \ldots, \alpha_k$. Let $\sigma(\langle \alpha \rangle)=\langle \beta \rangle$ and $\sigma(\langle \alpha_i \rangle)=\langle \beta_i \rangle$ for $1\leq i\leq k$. From Lemma 2.3 (iii), $\sigma(W)=\langle \beta_1, \ldots, \beta_k\rangle$. If $\sigma(\langle \alpha \rangle)$ is not a subspace of $\sigma(W)$, i.e., $\langle \beta \rangle$ is not a subspace of $\langle \beta_1, \ldots, \beta_k\rangle$, then $\beta, \beta_1, \ldots, \beta_k$ are linearly independent. Let $W'=\langle \beta, \beta_1, \ldots, \beta_k\rangle$ be a $(k+1)$-dimensional subspace. As $\sigma^{-1}$ is also an automorphism of $\mathcal{G}(\mathbb{V})$, applying $\sigma^{-1}$ to $W'$ together with Lemma 2.3, we obtain that
$$\sigma^{-1}(W')=\langle \alpha, \alpha_1, \ldots, \alpha_k \rangle=\langle \alpha_1, \ldots, \alpha_k \rangle=W,$$
which contradicts with Lemma 2.2.\hfill$\square$

By $I$ and $E_{ij}$, we respectively denote the identity matrix of order $n$ and the matrix unit of order $n$ with the $(i,j)$-entry 1 and others 0. We use $P_{kt}$ to denote the permutation matrix obtained from $I$ by permuting the $k$-th row and the $t$-th row (particularly, $P_{kk}=I$).
For an 1-dimensional subspace $\langle \alpha \rangle$, since $\langle \alpha \rangle=\langle a\alpha \rangle$ with $0\neq a\in F_q$, the expression can be not unique, but if we let the first nonzero component of $\alpha$ be 1, then the expression become unique. For convenience, in the following, we always suppose that each 1-dimensional subspace is of this form.

\begin{lemma}\ Let $\sigma$ be an automorphism of $\mathcal{G}(\mathbb{V})$. Then there exists an invertible matrix $B$ such that the automorphism $\sigma_B\circ \sigma$ fixes each $\langle e_i \rangle$ ($1\leq i\leq n$).
\end{lemma}

\noindent
{\bf Proof.} \ Let $\sigma(\langle e_1 \rangle)=\langle \alpha \rangle$ and $\alpha=\sum_{i=k}^{n}a_i e_i$ with $a_k=1$ and $1\leq k\leq n$. Then we take
$$B_1=P_{1k}(I-\sum_{j=k+1}^n a_jE_{jk}).$$
It is obvious that $\sigma_{B_1}\circ \sigma(\langle e_1 \rangle)=\langle e_1 \rangle$. Set $\sigma_1=\sigma_{B_1}\circ \sigma$. Further, suppose $\sigma_1(\langle e_2 \rangle)=\langle \beta \rangle$. Since $e_1$ and $e_2$ are linearly independent, so are $e_1$ and $\beta$ from Lemma 2.3 (ii). Then we can write $\beta=\sum_{i=1}^{n}b_i e_i$ and $\sum_{i=2}^{n}|b_i|\neq 0$. If $b_t\neq 0$ ($t\geq 2$), then we take $$B_2=P_{2t}(I-b_t^{-1}\sum_{i\neq t}b_i E_{it}),$$
and one can easily check that $\sigma_{B_2}\circ \sigma_1$ fixes $\langle e_2 \rangle$. Set $\sigma_2=\sigma_{B_2}\circ \sigma_1$. Note that $\sigma_2$ also fixes $\langle e_1 \rangle$.

Let $W=\langle e_1, e_2 \rangle$ and $\sigma_2(\langle e_3 \rangle)=\langle \gamma \rangle$.
It follows from Lemma 2.3 that $\sigma_2(W)=W$. Since $e_3\notin W$, by Lemma 2.4 (ii), we have $\langle \gamma \rangle$ is not a subspace of $W$ (i.e., $\gamma \notin W$). Thus we can write $\gamma=\sum_{i=1}^{n}c_i e_i$ and $\sum_{i=3}^{n}|c_i|\neq 0$. Similarly, we can take an invertible matrix $B_3$ such that $\sigma_{B_3}\circ \sigma_2$ fixes $\langle e_3 \rangle$ and $\langle e_1 \rangle$, $\langle e_2 \rangle$. Proceeding in this method, we can find a matrix sequence $B_1, B_2, \ldots, B_n$ such that
$\sigma_{B_n}\circ \sigma_{B_{n-1}}\circ \cdots \circ \sigma_{B_1}\circ \sigma$ fixes every $\langle e_i \rangle$ ($1\leq i\leq n$).
Let $B=B_n B_{n-1}\cdots B_1$, then $\sigma_{B}\circ \sigma$ is what we want.
\hfill$\square$

\section{ Automorphisms of the subspace sum graph $\mathcal{G}(\mathbb{V})$ }

\quad In this section, we characterize the automorphisms of the subspace sum graph $\mathcal{G}(\mathbb{V})$. Let $\mathbb{V}$ be a finite dimensional vector space over a finite base field $F_q$ with $dim(\mathbb{V})\geq 3$.

\begin{theorem} \ A mapping $\sigma$ on the vertex set $V$ of $\mathcal{G}(\mathbb{V})$ is an automorphism of $\mathcal{G}(\mathbb{V})$ if and only if $\sigma$ can be uniquely decomposed as $\sigma=\sigma_A\circ \sigma_f$, where $\sigma_A$ and $\sigma_f$ are two standard automorphisms defined as above.
\end{theorem}

\noindent
{\bf Proof.} \ The proof for the sufficiency part is clear. In what follows, we present the proof for the necessity part. From Lemma 2.5, we can find an invertible matrix $B$ such that $\sigma_B\circ \sigma$ fixes every $\langle e_i \rangle$ ($1\leq i\leq n$).
Let $\sigma_1=\sigma_B\circ \sigma$. The remaining proof is divided into the following claims.

\vskip 3pt
\noindent
{\sf Claim 1.} \ {\it For a nonzero vector $\alpha=\sum_{i=1}^{n}a_i e_i$ with $a_i\in F_q$, assume
$\sigma_1(\langle \alpha \rangle)=\langle \beta \rangle$, where $\beta=\sum_{i=1}^{n}b_i e_i$ with $b_i\in F_q$. Then $a_i=0$ if and only if $b_i=0$ ($1\leq i\leq n$).}

Let $U_j=\langle e_1, \ldots, e_{j-1}, e_{j+1}, \ldots, e_n \rangle$ for $1\leq j\leq n$. Since $\sigma_1$ fixes every
$\langle e_i \rangle$ for $1\leq i\leq n$, then from Lemma 2.3 (iii), $\sigma_1(U_j)=U_j$ for $1\leq j\leq n$.
If $a_j=0$, then $\langle \alpha \rangle\subset U_j$. Further, from Lemma 2.4 (i), $\langle \beta \rangle \subset U_j$, which implies that $b_j=0$. By applying $\sigma_1^{-1}$ to the above, we obtain $a_j=0$ whenever $b_j=0$.

\vskip 3pt
It follows from Claim 1 that $\sigma_1(\langle e_i+ae_j \rangle)=\langle e_i+be_j \rangle$ with $a, b\in F_q$ and $1\leq i<j\leq n$. Thus we define a function $f_{ij}$ over $F_q$ such that $\sigma_1(\langle e_i+ae_j \rangle)=\langle e_i+f_{ij}(a)e_j \rangle$ and  $f_{ij}(x)=0$ if and only if $x=0$. Now we investigate the properties of $f_{ij}$.

\vskip 3pt
\noindent
{\sf Claim 2.} \ {\it $\sigma_1(\langle e_i+\sum_{j=i+1}^{n}a_j e_j \rangle)=\langle e_i+\sum_{j=i+1}^{n}f_{ij}(a_j)e_j \rangle$ for $1\leq i<j\leq n$.}

From Claim 1, suppose $\sigma_1(\langle e_i+\sum_{j=i+1}^{n}a_j e_j \rangle)=\langle e_i+\sum_{j=i+1}^{n}b_j e_j \rangle$.
For $i<k\leq n$, let $$W_k=\langle e_i+a_ke_k, e_{i+1}, \ldots, e_{k-1}, e_{k+1}, \ldots, e_n \rangle.$$
Then from Lemmas 2.3 and 2.5,
$$\sigma_1(W_k)=\langle e_i+f_{ik}(a_k)e_k, e_{i+1}, \ldots, e_{k-1}, e_{k+1}, \ldots, e_n \rangle.$$
Since $\langle e_i+\sum_{j=i+1}^{n}a_j e_j \rangle\subset W_k$, then after applying $\sigma_1$, we derive from Lemma 2.4(i) that
$$\langle e_i+\sum_{j=i+1}^{n}b_j e_j \rangle \subset \langle e_i+f_{ik}(a_k)e_k, e_{i+1}, \ldots, e_{k-1}, e_{k+1}, \ldots, e_n \rangle,$$
which indicates that $b_k=f_{ik}(a_k)$. Notice that $i<k\leq n$, then the conclusion follows.

\vskip 3pt
\noindent
{\sf Claim 3.} \ {\it Suppose $2\leq i<j\leq n$ and $\forall a, b\in F_q$, then}
\begin{spacing}{1}
\begin{enumerate}
\item[(i)]  \textit{$f_{1j}(ab)=f_{1i}(a)f_{ij}(b)$};
\item[(ii)]  \textit{$f_{1j}(a)=f_{1i}(a)f_{ij}(1)=f_{1i}(1)f_{ij}(a)$};
\item[(iii)]  \textit{$\frac{f_{1j}(a)}{f_{1j}(1)}=\frac{f_{1i}(a)}{f_{1i}(1)}=\frac{f_{ij}(a)}{f_{ij}(1)}$}.
\end{enumerate}
\end{spacing}

Applying $\sigma_1$ on $\langle e_1 \rangle \subset \langle e_1+ae_i+abe_j, e_i+be_j \rangle$, together with Lemmas 2.4 and 2.5 and Claim 2, we obtain $$\langle e_1 \rangle \subset \langle e_1+f_{1i}(a)e_i+f_{1j}(ab)e_j, e_i+f_{ij}(b)e_j \rangle,$$
which points out that $f_{1j}(ab)=f_{1i}(a)f_{ij}(b)$. Taking respectively $a=1$, $b=1$ and $a=b=1$ in (i), one can easily deduce the results of (ii) and (iii).

\vskip 3pt
\noindent
{\sf Claim 4.} \ {\it  Let $f=\frac{f_{12}}{f_{12}(1)}$, then for $\forall a, b\in F_q$ and $0\neq c\in F_q$,}
\begin{spacing}{1}
\begin{enumerate}
\item[(i)]  \textit{$f(ab)=f(a)f(b)$};
\item[(ii)]  \textit{$f(1)=1$ and $f(c^{-1})=f(c)^{-1}$};
\item[(iii)]  \textit{$f(-1)=-1$ and $f(-a)=-f(a)$};
\item[(iv)]  \textit{$f(a+b)=f(a)+f(b)$};
\item[(v)]  \textit{$f$ is an automorphism of $F_q$}.
\end{enumerate}
\end{spacing}

From Claim 3 (iii), we see for $2\leq i<j\leq n$ that $$f(a)=\frac{f_{12}(a)}{f_{12}(1)}=\frac{f_{1j}(a)}{f_{1j}(1)}=\frac{f_{ij}(a)}{f_{ij}(1)}.$$
Then
$$f(ab)=\frac{f_{1n}(ab)}{f_{1n}(1)}=\frac{f_{12}(a)f_{2n}(b)}{f_{12}(1)f_{2n}(1)}=f(a)f(b),$$
which proves (i).

It follows from the definition of $f$ that $f(1)=1$. Then $1=f(1)=f(cc^{-1})=f(c)f(c^{-1})$, and thus $f(c^{-1})=f(c)^{-1}$, which establishes (ii).

Since $\langle e_1-ae_3 \rangle \subset \langle e_1-ae_2, e_2-e_3 \rangle$, after applying $\sigma_1$ together with Lemma 2.4 and Claim 2, it follows that $\langle e_1+f_{13}(-a)e_3 \rangle \subset \langle e_1+f_{12}(-a)e_2, e_2+f_{23}(-1)e_3 \rangle$. Thus we obtain
$$f_{13}(-a)=-f_{12}(-a)f_{23}(-1),$$
which, together with Claim 3 (i), imply that $f_{13}(-a)=-f_{13}(a)$.  Further, we obtain
$$f(-a)=-f(a),$$ and $f(-1)=-f(1)=-1$ if $a=1$. The proof of (iii) is finished.

Note that $f(0)=0$, then if $a=0$ in (iv), the result is clear. Now suppose $a\neq 0$. Applying $\sigma_1$ on $\langle e_1+(a+b)e_n \rangle \subset \langle e_1-ae_2+ae_n, e_2+a^{-1}be_n \rangle$, together with Lemma 2.4 and Claim 2, we deduce
$$\langle e_1+f(a+b)f_{1n}(1)e_n \rangle \subset \langle e_1+f(-a)f_{12}(1)e_2+f(a)f_{1n}(1)e_n, e_2+f(a^{-1}b)f_{2n}(1)e_n \rangle,$$
which indicates that
\begin{equation*}
\begin{array}{rcl}
   f(a+b)f_{1n}(1)&=& f(a)f_{1n}(1)-f(-a)f_{12}(1)f(a^{-1}b)f_{2n}(1)\\
   &=& f(a)f_{1n}(1)+f(a)f_{12}(1)f(a^{-1}b)f_{2n}(1)\\
   &=& f(a)f_{1n}(1)+[f(a)f(a^{-1}b)][f_{12}(1)f_{2n}(1)]\\
   &=& f(a)f_{1n}(1)+f(b)f_{1n}(1).
\end{array}
\end{equation*}
As a consequence, $f(a+b)=f(a)+f(b)$ holds.

Combining (i)-(iv), (v) follows immediately.

\vskip 3pt
\noindent
{\sf Claim 5.}\ {\it There exists a diagonal matrix $D$ such that $\sigma_f^{-1}\circ \sigma_D\circ \sigma_1$ fixes each 1-dimensional subspace.}

From Claims 3 and 4, $f_{ij}(a_{j})=f(a_{j})f_{ij}(1)=f(a_{j})\frac{f_{1j}(1)}{f_{1i}(1)}$. Then the result of Claim 2 can be rewritten as
\begin{equation*}
\sigma_1(\langle e_i+\sum_{j=i+1}^{n}a_j e_j \rangle)= \langle e_i+\sum_{j=i+1}^{n}f_{ij}(a_j)e_j \rangle
= \langle e_i+\frac{1}{f_{1i}(1)}\sum_{j=i+1}^{n}f(a_j)f_{1j}(1)e_j \rangle.
\end{equation*}
Let $f_{11}(1)=1$ and $D=diag(1, f_{12}^{-1}(1), \ldots, f_{1n}^{-1}(1))$, then
$$\sigma_D\circ \sigma_1(\langle e_i+\sum_{j=i+1}^{n}a_j e_j \rangle)=\langle e_i+\sum_{j=i+1}^{n}f(a_j) e_j \rangle,$$
and further it is clear that $\sigma_f^{-1}\circ \sigma_D\circ \sigma_1$ fixes each 1-dimensional subspace.
Set $\sigma_2=\sigma_f^{-1}\circ \sigma_D\circ \sigma_1$.

Now we are in a position to complete the proof.

Let $W\in V$ be any $k$-dimensional nontrivial proper subspace of $\mathbb{V}$ and $\{w_1, w_2, \ldots, w_k\}$ a basis of $W$ with $k\geq 2$. From Claim 5, $\sigma_2$ fixes every 1-dimensional subspace, then $\sigma_2$ fixes $\langle w_i \rangle$ ($1\leq i\leq k$). Using Lemma 2.3, we obtain $\sigma_2$ fixes $W$, i.e., $\sigma_2$ fixes every vertex of $\mathcal{G}(\mathbb{V})$. As a result, $\sigma_2=\sigma_f^{-1}\circ \sigma_D\circ \sigma_B\circ \sigma$  is the identity mapping. Thus
$$\sigma=\sigma_{B^{-1}}\circ \sigma_{D^{-1}} \circ \sigma_f=\sigma_A\circ \sigma_f$$ with $A=B^{-1}D^{-1}$.
Next, we show the decomposition of $\sigma$ is unique. Suppose there are two distinct decompositions,
$\sigma_{A_1}\circ \sigma_{f_1}=\sigma_{A_2}\circ \sigma_{f_2}$, which implies that $\sigma_{A_2^{-1}A_1}=\sigma_{f_2f_1^{-1}}$.
Since $\sigma_{f_2f_1^{-1}}$ fixes each $\langle e_i \rangle$, so does $\sigma_{A_2^{-1}A_1}$. Then we see $A_2^{-1}A_1$ is a diagonal matrix. Similarly, $\sigma_{A_2^{-1}A_1}$ also fixes every $\langle e_i+e_j \rangle$ for $1\leq i<j\leq n$, from which it follows that $A_2^{-1}A_1$ is a scalar matrix. So, $A_2$ is a nonzero scalar multiple of $A_1$, and thus $\sigma_{A_1}=\sigma_{A_2}$. As a consequence, $\sigma_{f_1}=\sigma_{f_2}$.

The proof of Theorem 3.1 is completed.\hfill$\square$

\vskip 4pt
Note that the automorphism group of $F_q$ with $q=p^m$ is a cyclic group of order $m$. Then the following result follows from Theorem 3.1 immediately.

\begin{corol} \ The automorphism group of $\mathcal{G}(\mathbb{V})$ is isomorphic to $PGL_n(F_q)\times \mathbb{Z}_m$, where $PGL_n(F_q)$ is the quotient group of all the $n\times n$ invertible matrices over $F_q$ to the normal subgroup of all the nonzero scalar matrices over $F_q$.
\end{corol}

\vskip 5pt
\noindent
{\bf Acknowledgement}

The authors thank the anonymous referees for a careful and thorough reading of the paper and for valuable comments, and also thank the support of the Fundamental Research Funds for the Central Universities (No. 2017BSCXB53).

\end{document}